\documentclass[12pt]{article}
\usepackage{amsmath,amsthm,amssymb}
\begin{document}
\newtheorem{thm}{Theorem}
\newtheorem{pro}[thm]{Proposition}
\newtheorem{cor}[thm]{Corollary}
\newtheorem{lem}[thm]{Lemma}
\newtheorem{df}[thm]{Definition}
\newtheorem{rem}[thm]{Remark}
\newtheorem{conj}[thm]{Conjecture}
\newcommand{\mset}{\mathop{\mathrm{maxset\hspace{1mm}}}\nolimits}
\newcommand{\lset}{\mathop{\mathrm{minset\hspace{1mm}}}\nolimits}
\newcommand{\Leg}{\mathop{\mathrm{Length}}\nolimits}
\title{Geodesics of Hofer's metric on the space of Lagrangian
submanifolds}
\author{Hiroshi Iriyeh and Takashi Otofuji}
\date{}
\maketitle
\begin{abstract}
We study geodesics of Hofer's metric on the space of Lagrangian submanifolds
in arbitrary symplectic manifolds from the variational point of view.
We give a characterization of length-critical paths with respect to this
metric.
As a result, we prove that if two Lagrangian submanifolds are disjoint
then we cannot join them by length-minimizing geodesics.
\end{abstract}

\section{Introduction}

\indent\indent
The study of geodesics on Hamiltonian diffeomorphism groups with respect
to Hofer's metric was initiated by Bialy and Polterovich \cite{BP}
in the case of the standard symplectic space $\mathbb{R}^{2n}$
and many deep results are known today for general symplectic manifolds
(see, for example, \cite{LM,LM2,Os,Us,Po}).

By contrast there are few results on Hofer's geodesics on the space of
Lagrangian submanifolds which are Hamiltonian isotopic to one.
Milinkovi\'c \cite{Mi,Mi2} studied the case of Hamiltonian isotopy class of
the zero section of the cotangent bundle of a compact manifold.
Akveld-Salamon \cite{AS} proved length minimizing properties of some
Lagrangian loops of $\mathbb{R}P^n$ in $\mathbb{C}P^n$.

As for a related topic,
recently Ostrover \cite{Os} compared the two Hofer geometries above.
By corresponding a Hamiltonian diffeomorphism of a closed symplectic manifold
$(M,\omega)$ to its graph, one can embed the Hamiltonian diffeomorphism group
$\mathrm{Ham}(M,\omega)$ into the space $\mathcal{L}$ of all Lagrangian
submanifolds of $(M \times M, -\omega \oplus \omega)$ which are Hamiltonian
isotopic to the diagonal.
He showed that the image of $\mathrm{Ham}(M,\omega)$ in $\mathcal{L}$ is
strongly distorted.

In the present paper, we shall initiate the study of geodesics on the space
of Lagrangian submanifolds in arbitrary symplectic manifolds
from the variational point of view.

Let $(M,\omega)$ be a connected symplectic manifold without boundary.
Let $L \subset M$ be a closed Lagrangian submanifold.
A smooth family $\{ L_t \}_{0 \leqq t \leqq 1}$ of Lagrangian submanifolds
of $M$, such that $L_t$ is diffeomorphic to $L$,
is called an {\it exact (Lagrangian) path} connecting $L_0$ and $L_1$,
if there exists a Hamiltonian isotopy $\psi_t$ of $M$ such that
$\psi_{t}(L_0)=L_t$ for every $t \in [0,1]$.
The Hofer length of an exact path $\{ L_t \}$ is given by
 \begin{eqnarray*}
  \Leg(\{ L_t \})=
   \int_{0}^{1}\left( \max_{p \in L_t}H(t,p)-\min_{p \in L_t}H(t,p) \right)dt,
 \end{eqnarray*}
where $H \in C_{0}^{\infty}([0,1] \times M)$ is the compactly supported
Hamiltonian function which generates a Hamiltonian isotopy
$\psi_t: M \to M$ satisfying $\psi_{t}(L_0)=L_t$
(see Section 2 for the details).

We denote by $\mathcal{L}=\mathcal{L}(M,\omega,L)$ the space of all
Lagrangian submanifolds of $M$ which are connected to $L$ by an exact path.
For any $L_0,L_1 \in \mathcal{L}(M,\omega,L)$, a bi-invariant
pseudodistance $d$ on $\mathcal{L}$ is defined by
 \begin{eqnarray*}
  d(L_0,L_1)=\inf \Leg(\{ L_t \}),
 \end{eqnarray*}
where the infimum is taken over all exact paths on $\mathcal{L}$ connecting
$L_0$ and $L_1$.

Chekanov \cite{Ch1,Ch2} proved that $d$ is non-degenerate
if $(M,\omega)$ is {\it tame}, i.e.,
it admits an almost complex structure $J$ such that the bilinear form
$\omega(\cdot,J\cdot)$ defines a complete Riemannian metric on $M$
with bounded curvature and with injectivity radius bounded away from zero.
The distance $d$ is called {\it Hofer's distance} on $\mathcal{L}$.

An exact path $\{ L_t \}$ is said to be {\it l-critical} if it is a critical
point of the Hofer length functional
(see Section 3 for the precise definition).
Our main result is the following.
\begin{thm}
Let $(M,\omega)$ be a symplectic manifold without boundary and $L$ a
closed Lagrangian submanifold of $M$.
Let $\{ L_t \}_{0 \leqq t \leqq 1} \subset \mathcal{L}(M,\omega,L)$
be the exact path defined by $L_t=\psi_{t}(L_0)$, where $\psi_t$ is
the Hamiltonian isotopy which is generated by a Hamiltonian function
$H \in C_{0}^{\infty}([0,1] \times M)$.
Then, $\{ L_t \}$ is {\it l-critical} if and only if there exist two points
$p_{+},p_{-} \in \bigcap_{t \in [0,1]} L_t \subset M$
satisfying
 \begin{eqnarray*}
  \max \left(H|_{L_t} \right)=H(t,p_+),\quad
  \min \left(H|_{L_t} \right)=H(t,p_-)
 \end{eqnarray*}
for all $t \in [0,1]$. \hfill \qed
\end{thm}
This result is a generalization of the result by Lalonde-McDuff in the case
of paths in a Hamiltonian diffeomorphism groups \cite[Theorem 1.9]{LM}.

\begin{thm}[Lalonde-McDuff]
Let $(M,\omega)$ be a symplectic manifold without boundary.
Let $\{ \psi_t \}_{0 \leqq t \leqq 1}$
be the Hamiltonian isotopy which is generated by a Hamiltonian function
$H \in C_{0}^{\infty}([0,1] \times M)$.
Then, $\{ \psi_t \}$ is a critical point of the Hofer length functional
if and only if there exist two points $p_{+},p_{-} \in M$
satisfying
 \begin{eqnarray*}
  \max_{p \in M} H(t,p)=H(t,p_+),\quad
  \min_{p \in M} H(t,p)=H(t,p_-)
 \end{eqnarray*}
for all $t \in [0,1]$. \hfill \qed
\end{thm}

This paper is organized as follows.
In Section 2 we review some facts about the Hofer geometry on the space
of Lagrangian submanifolds.
In Section 3 we introduce the notion of exact variations of an exact
Lagrangian path and describe the length function for such a variation.
The main theorem above is proved in Section 4.
In the last section we give some applications of our main theorem.
In particular, we prove that if two Lagrangian submanifolds are disjoint
then we cannot join them by length-minimizing geodesics.

\section{The Hofer length for exact Lagrangian paths}

\indent\indent
In this section we give a precise description of the Hofer geometry on the
space of Lagrangian submanifolds along Akveld-Salamon's treatment in
\cite[Section 2]{AS}.

Let $(M,\omega)$ be a connected symplectic manifold without boundary.
Let $L \subset M$ be a closed Lagrangian submanifold.
We denote by
 \begin{eqnarray*}
  \widehat{\mathrm{Lag}}(L,M)
   =\{ \iota \in \mathrm{Emb}(L,M)\ |\ \iota^{*}\omega=0 \}
 \end{eqnarray*}
the space of Lagrangian embeddings of $L$ into $M$. 
The group of diffeomorphisms $\mathrm{Diff}(L)$ of $L$ acts on
$\widehat{\mathrm{Lag}}(L,M)$ by $\iota \mapsto \iota \circ \phi$ for
$\phi \in \mathrm{Diff}(L)$.
Two Lagrangian embeddings $\iota_0, \iota_1 \in \widehat{\mathrm{Lag}}(L,M)$
belong to the same $\mathrm{Diff}(L)$-orbit if and only if they have the same
image.
Therefore, the quotient space
$\mathrm{Lag}(L,M):=\widehat{\mathrm{Lag}}(L,M)/\mathrm{Diff}(L)$
can be considered as the space of Lagrangian submanifolds of $M$ which are
diffeomorphic to $L$.
A path $\{ L_t \}_{0 \leqq t \leqq 1} \subset \mathrm{Lag}(L,M)$ is said to
be {\it smooth} if there exists a smooth map
$[0,1] \times L \to M: (t,x) \mapsto \iota_t(x)$ such that
$\iota_{t}(L)=L_t$ for all $t \in [0,1]$.
This $\{ \iota_t \}_{0 \leqq t \leqq 1}$ is called a {\it lift} of $\{ L_t \}$.

Next we explain that one can think of $\mathrm{Lag}(L,M)$ as an infinite
dimensional manifold.
Let $[0,1] \times L \to M: (t,x) \mapsto \iota_t(x)$ be a smooth map
such that $\iota_t \in \widehat{\mathrm{Lag}}(L,M)$ for all $t \in [0,1]$.
Let us introduce a $1$-form on $L$ defined by
 \begin{eqnarray}
  \alpha_t:=\omega\Bigl(\frac{d}{dt}\iota_{t},d\iota_{t} \cdot \Bigl)\ 
   \in \Omega^{1}(L).
 \end{eqnarray}
Then we have
 \begin{eqnarray*}
  0 = \frac{\partial}{\partial t}(\iota_{t}^{*}\omega) = d \alpha_{t}
 \end{eqnarray*}
and hence the tangent space of $\widehat{\mathrm{Lag}}(L,M)$ at $\iota$ is
given by
 \begin{eqnarray*}
  T_{\iota} \widehat{\mathrm{Lag}}(L,M)
   =\{ v \in C^{\infty}(L,\iota^{*}TM)\ |\ \omega(v,d\iota \cdot)
    \in \Omega^{1}(L): \mathrm{closed} \}.
 \end{eqnarray*}
The tangent space to the $\mathrm{Diff}(L)$-orbit at $\iota$ is described as
 \begin{eqnarray*}
  T_{\iota} (\iota \cdot \mathrm{Diff}(L))
   =\{ v=d\iota \circ \xi\ |\ \xi \in \mathcal{X}(L) \},
 \end{eqnarray*}
where $\mathcal{X}(L)$ denotes the space of vector fields on $L$.
The map $v \mapsto \omega(v,d\iota \cdot)$ induces the linear map
 \begin{eqnarray*}
  T_{\iota} \widehat{\mathrm{Lag}}(L,M)/
   T_{\iota} (\iota \cdot \mathrm{Diff}(L))
    \longrightarrow \{ \beta \in \Omega^{1}(L)\ |\ d\beta=0 \},
 \end{eqnarray*}
which is an isomorphism since $\iota: L \to M$ is Lagrangian.

Take $L' \in \mathrm{Lag}(L,M)$ and
two paths $\{ \iota_t \}$ and $\{ \iota'_t \}$ in
$\widehat{\mathrm{Lag}}(L,M)$ which satisfy $\iota'_t=\iota_t \circ \phi_t$
for some path $\{ \phi_t \} \subset \mathrm{Diff}(L)$ and
$[\iota_0]=[\iota'_0]=L'$.
The derivative of this condition yields
 \begin{eqnarray*}
  \frac{d}{dt}\iota'_t
   =\frac{d}{dt}\iota_t \circ \phi_t + d\iota_t \circ \xi_t \circ \phi_t,
 \end{eqnarray*}
where $\xi_t \in \mathcal{X}(L)$ generates $\phi_t \in \mathrm{Diff}(L)$
via $\frac{d}{dt}\phi_t = \xi_t \circ \phi_t$.
Using the fact that $\iota_t: L \to M$ is Lagrangian, the $1$-forms
$\alpha_t:=\omega(\frac{d}{dt}\iota_{t},d\iota_{t} \cdot)$ and
$\alpha'_t:=\omega(\frac{d}{dt}\iota'_{t},d\iota'_{t} \cdot)$
satisfy
 \begin{eqnarray*}
  \alpha'_t = \phi_{t}^{*}\alpha_t.
 \end{eqnarray*}
Hence closed $1$-forms $\alpha, \alpha' \in \Omega^{1}(L)$ associated
with Lagrangian embeddings $\iota$ and $\iota'=\iota \circ \phi$ represent
the same tangent vector of $\mathrm{Lag}(L,M)$ at $L'$ if and only if
$\alpha' = \phi^{*}\alpha$ (or, $\iota_{*}\alpha=\iota'_{*}\alpha'$).

Summarizing these arguments, we have
\begin{lem}
The tangent space of $\mathrm{Lag}(L,M)$ at $L'$ can be identified with
the space of closed $1$-forms on $L$, i.e.,
 \begin{eqnarray*}
  T_{L'} \mathrm{Lag}(L,M) = \{ \beta \in \Omega^{1}(L)\ |\ d\beta=0 \}.
 \end{eqnarray*}
\end{lem}

Let $[0,1] \to \mathrm{Lag}(L,M): t \mapsto L_t$ be a smooth path of
Lagrangian submanifolds.
We define the velocity vector of the path $\{ L_t \}$ at time $t$ by
 \begin{eqnarray*}
  \frac{d}{dt}L_t := {\iota_t}_{*} \alpha_t,
 \end{eqnarray*}
where $\{ \iota_t \}$ is a lift of $\{ L_t \}$ and $\alpha_t$ is the $1$-form
defined by equation $(1)$.
The above argument proves that
$\beta_t:={\iota_t}_{*} \alpha_t \in \Omega^{1}(L_t)$
is closed and independent of the choice of the lift $\{ \iota_t \}$.

Our main target of this research is Hamiltonian isotopies of Lagrangian
submanifolds.
As we shall see below, a smooth path $\{ L_t \}_{0 \leqq t \leqq 1}$ in
$\mathrm{Lag}(L,M)$ whose velocity vectors $\beta_t \in \Omega^{1}(L_t)$
are exact for all $t \in [0,1]$ is nothing but an exact path.
Remember that $\mathcal{L}=\mathcal{L}(M,\omega,L)$ consists of all
Lagrangian submanifolds of $M$ connected to $L$ by an exact path,
hence $\mathcal{L}$ is a subset of $\mathrm{Lag}(L,M)$.
Let $\{ \psi_t \}_{0 \leqq t \leqq 1}$ be the Hamiltonian isotopy which is
generated by a Hamiltonian function $H \in C_{0}^{\infty}([0,1] \times M)$,
i.e.,
 \begin{eqnarray*}
  \frac{d}{dt}\psi_t = X_t \circ \psi_t,\ \omega(X_t,\cdot)=dH_t,\
   \psi_0=id.
 \end{eqnarray*}
\begin{lem}[\cite{AS}]
Let $\{ L_t \}_{0 \leqq t \leqq 1} \subset \mathrm{Lag}(L,M)$ be a smooth
path of Lagrangian submanifolds and $\{ \psi_t \}$ be the Hamiltonian isotopy
of $M$ generated by a Hamiltonian function $H$.
Then $L_t=\psi_{t}(L_0)$ for every $t \in [0,1]$ if and only if
 \begin{eqnarray*}
  \frac{d}{dt}L_t = d \left( H_{t}|_{L_t} \right)
 \end{eqnarray*}
for every $t \in [0,1]$. \hfill \qed
\end{lem}
Moreover, we have

\begin{lem}[\cite{AS}]
Let $\{ L_t \}_{0 \leqq t \leqq 1} \subset \mathrm{Lag}(L,M)$ be a smooth
path of Lagrangian submanifolds.
Then $\frac{d}{dt}L_t \in \Omega^{1}(L_t)$ is exact for every $t \in [0,1]$
if and only if $\{ L_t \}$ is an exact path, that is,
there exists a Hamiltonian isotopy $\{ \psi_t \}$
such that $\psi_{t}(L_0)=L_t$ for every $t \in [0,1]$. \hfill \qed
\end{lem}
We refer the reader to \cite[Lemmas 2.2 and 2.3]{AS} for proofs of
Lemmas 4 and 5.
Lemma 5 reads that a smooth path
$\{ L_t \}_{0 \leqq t \leqq 1} \subset \mathrm{Lag}(L,M)$
is an exact path if and only if there exists a lift
$\{ \iota_t \}_{0 \leqq t \leqq 1}$ of $\{ L_t \}$ which satisfies
 \begin{eqnarray}
  \alpha_{t}:=\omega\Bigl(\frac{d}{dt}\iota_{t},d\iota_{t} \cdot \Bigl)=dh_{t}
 \end{eqnarray}
for a smooth function $h \in C^{\infty}([0,1] \times L)$.
We call this smooth function $h$ an {\it associated function} of
$\{ \iota_t \}$.

Any $\{ \psi_t \}$ satisfying the condition in Lemma 5 is generated by
a Hamiltonian function $H \in C_{0}^{\infty}([0,1] \times M)$ such that
$H_t \circ \iota_t = h_t$ for each $t \in [0,1]$ and vice versa.

For any Lagrangian submanifold $L$, we can define the length of smooth paths
on the Hamiltonian isotopy class $\mathcal{L}=\mathcal{L}(M,\omega,L)$.
The {\it Hofer length} of an exact path
$\{ L_t \}_{0 \leqq t \leqq 1} \subset \mathcal{L}$ is defined by
 \begin{eqnarray*}
  \Leg(\{ L_t \}):=\int_{0}^{1} \Bigl\lVert \frac{d}{dt}L_t \Bigl\rVert dt   
  =\int_{0}^{1}\left( \max_{x \in L}h(t,x) - \min_{x \in L}h(t,x) \right)dt.
 \end{eqnarray*}
This quantity is independent of the choice of the lift $\{ \iota_t \}$ of
$\{ L_t \}$ and its associated function $h$.
In particular, Akveld-Salamon \cite{AS} constructed a Hamiltonian function
$H$ which satisfies
 \begin{eqnarray}
  \max_{M}H_{t}=\max_{L}h_{t},\quad \min_{M}H_{t}=\min_{L}h_{t}
 \end{eqnarray}
for each $t$.
Later, we must use more refined construction of Hamiltonian functions
for a proof of our main theorem (see Section 4).
By the above arguments and equation $(3)$, we have
 \begin{eqnarray*}
  \Leg(\{ L_t \})=\inf_{\psi_{t}(L_0)=L_t} l(\{ \psi_t \}),
 \end{eqnarray*}
where the infimum is taken over all Hamiltonian isotopies
$\{ \psi_t \}$ which satisfy
$\psi_{t}(L_0)=L_t$ for all $t \in [0,1]$ and $l(\{ \psi_t \})$ denotes
the Hofer length of the Hamiltonian path $\{ \psi_t \}$ (cf.\ \cite{LM,Po}).

\section{Length-critical paths}

\indent\indent
In this section we introduce the notion of exact variations of an exact path
and describe the length function for such a variation.

Let $\{ L_t \}_{0 \leqq t \leqq 1} \subset \mathcal{L}$ be a smooth
{\it regular} (i.e., $\frac{d}{dt}L_t \neq 0$ for all $t \in [0,1]$)
exact path connecting $L_0$ and $L_1$ and $I$ be an open interval in
$\mathbb{R}$ containing $0$.

\vspace{1mm}

\noindent
{\bf Definition.}\ A smooth map
$\Phi:I \times [0,1] \to \widehat{\mathrm{Lag}}(L,M)$
is called an {\it exact (Lagrangian) variation} of $\{ L_t \}$
if there exists a lift $\{ \iota_t \}_{0 \leqq t \leqq 1}$ of $\{ L_t \}$
which has the following properties:

$(\mathrm{i})$ $\Phi(0,t)=\iota_t$ for all $t \in [0,1]$;

$(\mathrm{ii})$ $\Phi(s,0)=\iota_0, \Phi(s,1)=\iota_1$ for all $s \in I$;

$(\mathrm{iii})$ $\Phi(s,t)=:\iota_{s,t}$ is a lift of an exact path for each
$s \in I$.

\vspace{1mm}

\noindent
For every $s \in I$ we denote by $l(s)$ the length of the exact path
$\{ \iota_{s,t}(L) \}_{0 \leqq t \leqq 1} \subset \mathcal{L}$:
 \begin{eqnarray*}
  l(s)=\Leg(\{ \iota_{s,t}(L) \}_{0 \leqq t \leqq 1}).
 \end{eqnarray*}

We call a function arising in this way a {\it length function}.
Since the length function $l(s)$ is not necessarily differentiable,
we must modify the notion of critical points of $l(s)$.
Our treatment here is motivated by Polterovich's book \cite[Chapter 12]{Po}.
\begin{pro}
For any length function $l(s)$, there exists a convex function $u(s)$
and positive numbers $\delta$, $C$ such that
 \begin{eqnarray}
  |l(s)-u(s)| \leqq Cs^2
 \end{eqnarray}
for all $s \in (-\delta,\delta)$.
\hfill \qed
\end{pro}
We shall give its proof later.
From this fact, we arrive at the following definition.

\vspace{1mm}

\noindent
{\bf Definition.}\ We say that the length function $l(s)$ has a
{\it critical point} at $s=0$ if it is a minimal point for a convex
function $u(s)$ which satisfies $(4)$.

\vspace{1mm}

\noindent
{\bf Remark.}\ This definition does not depend on the particular choice of
the convex function $u$ satisfying $(4)$.

\vspace{1mm}

Using this modified notion of critical points, we can define the notion of
$l$-critical paths.

\vspace{1mm}

\noindent
{\bf Definition.}\ A regular exact path
$\{ L_t \}_{0 \leqq t \leqq 1} \subset \mathcal{L}(M,\omega,L)$
is said to be {\it l-critical} if for every exact variation of $\{ L_t \}$
the length function $l(s)$ has a critical point at $s=0$.

\vspace{1mm}

Here we review some facts about Hamiltonian loops.
Consider the space $\mathcal{H}_0$ of all normalized Hamiltonians
$F: [0,1] \times M \to \mathbb{R}$ endowed with the Hofer norm
 \begin{eqnarray*}
  \lVert F \rVert
  =\int_{0}^{1}\left( \max_{x \in M}F(t,x) - \min_{x \in M}F(t,x) \right)dt.
 \end{eqnarray*}
Denote by $\{ \varphi_{t}^{F} \}_{0 \leqq t \leqq 1}$ the Hamiltonian flow
generated by $F$ with $\varphi_{0}^{F}=id$.
The Hamiltonian diffeomorphism group $\mathrm{Ham}^{c}(M,\omega)$ of $M$
is defined by
 \begin{eqnarray*}
  \mathrm{Ham}^{c}(M,\omega):=\{ \varphi_{1}^{F} |\ F \in \mathcal{H}_0 \}.
 \end{eqnarray*}

Next we describe exact variations of an exact path
$\{ L_t \}_{0 \leqq t \leqq 1}$.
Take a lift $\{ \iota_t \}_{0 \leqq t \leqq 1}$ of $\{ L_t \}$ and its
associated function $\{ h_t \}$.
Let $\{ \iota_{s,t} \}_{s \in I}$ be an exact variation of $\{ L_t \}$ with
their associated functions $\{ h_{s,t} \}_{s \in I}$.
We denote by $\{ L_{s,t} \}$ the exact paths $\{ \iota_{s,t}(L) \}$.
Since $\mathrm{Ham}^{c}(M,\omega)$ acts transitively on
$\mathcal{L}(M,\omega,L)$,
we can take Hamiltonian paths $\{ \psi_{s,t} \}$ and
$\{ \psi_t \} \subset \mathrm{Ham}^{c}(M,\omega)$
generated by Hamiltonian functions $\{ H_{s,t} \}$ and $\{ H_t \}$,
respectively, satisfying
 \begin{eqnarray}
  L_t = \psi_{t}(L_0),\quad L_{s,t} = \psi_{s,t}(L_0)
 \end{eqnarray}
and $\psi_{s,0} = id,\ \psi_0 = id,\ \psi_{0,t} = \psi_t$.
Define $\{ f_{s,t} \} \subset \mathrm{Ham}^{c}(M,\omega)$ by the equation
 \begin{eqnarray}
  \psi_{s,t} = f_{s,t}^{-1} \circ \psi_t.
 \end{eqnarray}

The following is the key lemma of a proof of our main theorem.
\begin{lem}
We can assume that the Hamiltonian path $\{ f_{s,t} \}_{0 \leqq t \leqq 1}$
is a loop for each parameter $s \in I$.
\end{lem}
{\it Proof.} The Hamiltonian path $\{ f_{s,t} \}_{0 \leqq t \leqq 1}$
starts from $id$ and satisfies $f_{s,1}(L_1)=L_1$ for each $s \in I$.
Let us introduce a new Hamiltonian path $\{ \psi'_{s,t} \}$ defined by
 \begin{eqnarray*}
  \psi'_{s,t}=\left( (\psi_{s,1} \circ \psi_{s,t}^{-1})^{-1} \circ f_{ts,1}
   \circ (\psi_{s,1} \circ \psi_{s,t}^{-1}) \right) \circ \psi_{s,t}.
 \end{eqnarray*}
Then we have
 \begin{eqnarray*}
  \psi'_{s,1} &=& ( id^{-1} \circ f_{s,1} \circ id) \circ \psi_{s,1}=\psi_1,\\
  \psi'_{s,0} &=& ( \psi_{s,1}^{-1} \circ f_{0,1} \circ \psi_{s,1})
   \circ \psi_{s,0}=\psi_{s,0}=id.
 \end{eqnarray*}
Here we define $\{ f'_{s,t} \} \subset \mathrm{Ham}^{c}(M,\omega)$ by the
equation $\psi'_{s,t} = (f'_{s,t})^{-1} \circ \psi_t$.
Then
 \begin{eqnarray*}
  f'_{s,1} &=& \psi_1 \circ (\psi'_{s,1})^{-1} =id,\\
  f'_{s,0} &=& \psi_0 \circ (\psi'_{s,0})^{-1} =id.
 \end{eqnarray*}
Hence $\{ f'_{s,t} \}$ is a Hamiltonian loop for every $s \in I$.

It suffices to show that $\psi'_{s,t}(L_0)=L_{s,t}$.
Indeed, we have
 \begin{eqnarray*}
  \psi'_{s,t}(L_0)
   &=& \left( (\psi_{s,1} \circ \psi_{s,t}^{-1})^{-1}
    \circ f_{ts,1} \circ (\psi_{s,1} \circ \psi_{s,t}^{-1}) \right)
     \circ \psi_{s,t}(L_0)\\
   &=& (\psi_{s,1} \circ \psi_{s,t}^{-1})^{-1}
    \circ f_{ts,1} \circ \psi_{s,1}(L_0)\\
   &=& (\psi_{s,1} \circ \psi_{s,t}^{-1})^{-1} \circ f_{ts,1}(L_1)\\
   &=& \psi_{s,t} \left( \psi_{s,1}^{-1}(L_1) \right) = L_{s,t}.
 \end{eqnarray*}
\hfill \qed

From now on, by virtue of Lemma 7, we assume that Hamiltonian paths
$\{ f_{s,t} \}_{0 \leqq t \leqq 1}$ defined by equation $(6)$ are loops.
Denote by $\mathcal{H}_1$ the subset of $\mathcal{H}_0$ which consists of
all Hamiltonians generating loops in $\mathrm{Ham}^{c}(M,\omega)$.
Denote by $\mathcal{V}$ the set of all smooth $1$-parameter families
$F(s,t,x)$ of functions from $\mathcal{H}_1$ with $F(0,t,x) \equiv 0$.
\begin{lem}
For an exact path $\{ L_t \}$, 
the set of length functions $l(s)$ associated with the exact variations of
$\{ L_t \}$ consists of all functions of the form
 \begin{eqnarray*}
  ||h-F(s)\circ\iota_t||,
 \end{eqnarray*}
where $F \in \mathcal{V}$ and $\{ \iota_t \}$ is any fixed lift of $\{ L_t \}$
with its associated function $h$.
Here, $||h||$ denotes the norm of the associated function
$h \in C^{\infty}([0,1] \times L)$ defined by
 \begin{eqnarray*}
  ||h||
  =\int_{0}^{1}\left( \max_{x \in L}h(t,x) - \min_{x \in L}h(t,x) \right)dt.
 \end{eqnarray*}
\end{lem}
{\it Proof.} For any exact variation $\{ \iota_{s,t} \}_{s \in I}$ of
$\{ L_t \}$ with their associated functions $\{ h_{s,t} \}_{s \in I}$,
by equation $(5)$, there exist paths
$\{ \psi_{s,t} \}_{s \in I} \subset \mathrm{Ham}^{c}(M,\omega)$ and
$\{ \phi_{s,t} \}_{s \in I} \subset \mathrm{Diff}(L)$ such that
 \begin{eqnarray}
  \psi_{s,t} \circ \iota_0 = \iota_{s,t} \circ \phi_{s,t}
 \end{eqnarray}
and
 \begin{eqnarray}
  \psi_t \circ \iota_0 = \iota_t \circ \phi_t \ \ \ (\mathrm{when}\ \ s=0).
 \end{eqnarray}
By equation $(6)$, we have
 \begin{eqnarray}
  H_{s,t}(y) = -F(s,t,f_{s,t}(y))+H(t,f_{s,t}(y))
 \end{eqnarray}
for all $y \in M$, where $\{ F_{s,t} \}$ are the family of Hamiltonian
functions generating the Hamiltonian isotopies $\{ f_{s,t} \}$ with
isotopy parameter $t \in [0,1]$.

Substituting $\iota_{s,t}(x),\ x \in L$, for $y$ in $(9)$,
we see by equations $(6),(7)$ and $(8)$ that
the associated functions $\{ h_{s,t} \}$ of $\{ \iota_{s,t} \}$ are
represented as
 \begin{eqnarray*}
  h_{s,t}(x) &=& H_{s,t} \circ \iota_{s,t}(x)\\
   &=& H(t,f_{s,t} \circ \iota_{s,t}(x))-F(s,t,f_{s,t} \circ \iota_{s,t}(x))\\
   &=& H(t,\iota_t \circ \phi_t \circ \phi_{s,t}^{-1}(x))
    -F(s,t,\iota_t \circ \phi_t \circ \phi_{s,t}^{-1}(x))\\
   &=& h_t(\phi_t \circ \phi_{s,t}^{-1}(x))
    -F(s) \circ \iota_t(\phi_t \circ \phi_{s,t}^{-1}(x)).
 \end{eqnarray*}
Since $\phi_t \circ \phi_{s,t}^{-1} \in \mathrm{Diff}(L)$, we have
 \begin{eqnarray*}
  l(s)=\Leg(\{ L_{s,t} \}_{0 \leqq t \leqq 1})=||h_{s}||
   =||h-F(s) \circ \iota_{t}||.
 \end{eqnarray*}

Conversely, if $\{ F(s) \}_s \in \mathcal{V}$ is given, then we can define
Hamiltonian functions $\{ H_{s,t} \}_s$ by equation $(9)$.
By integrating, we obtain Hamiltonian paths 
$\{ \psi_{s,t} \}_s \subset \mathrm{Ham}^{c}(M,\omega)$
satisfying equation $(6)$.
Consequently, the paths 
$\{ \psi_{s,t} \circ \iota_0 \}_s$ in $\widehat{\mathrm{Lag}}(L,M)$
define an exact variation whose length function is given by
$||h-F(s)\circ\iota_t||$. \hfill \qed

{\it Proof of Proposition 6.} By Lemma 8, any length function $l(s)$
associated with an exact variation of the exact path $\{ L_t \}$
is of the form
 \begin{eqnarray*}
  l(s) = ||h-F(s)\circ\iota_t||,
 \end{eqnarray*}
where $F \in \mathcal{V}$ and $\{ \iota_t \}$ is a lift of $\{ L_t \}$
with its associated function $h$.

Let us define a convex function $u(s)$ by
 \begin{eqnarray*}
  u(s) = ||h-sF'(0)\circ\iota_t||.
 \end{eqnarray*}
By the triangle inequality, we have
 \begin{eqnarray}
  l(s) &\leqq& ||h-sF'(0)\circ\iota_t||
   + ||sF'(0)\circ\iota_t-F(s)\circ\iota_t||\nonumber\\
       &=& u(s) + ||sF'(0)\circ\iota_t-F(s)\circ\iota_t||.
 \end{eqnarray}
Using Taylor's formula the inside of the second term of $(10)$ satisfies
 \begin{eqnarray*}
  sF'(0)\circ\iota_t-F(s)\circ\iota_t
   = -\frac{1}{2!} s^2 F''(c,t,\iota_t \cdot)
 \end{eqnarray*}
for all $s \in (-\delta,\delta)$,
where $c$ is a real constant between $0$ and $s$.
Hence we have
 \begin{eqnarray*}
  ||sF'(0)\circ\iota_t-F(s)\circ\iota_t||
   \leqq s^2
    \max_{x \in L,0 \leqq t \leqq 1, |s'| \leqq \delta}|F''(s',t,\iota_{t}(x))|
 \end{eqnarray*}
Setting $C=\max|F''|$ we get
 \begin{eqnarray*}
  l(s)-u(s) \leqq Cs^2.
 \end{eqnarray*}
Similarly, we have $u(s)-l(s) \leqq Cs^2$ with the same constant $C$.
\hfill \qed

\section{A Variational definition of geodesics}

\indent\indent
In this section we prove the main theorem stated in Section 1.
Let $\{ L_t \}_{0 \leqq t \leqq 1} \subset \mathcal{L}$ be a regular
exact path connecting $L_0$ and $L_1$.
The following notion is motivated by the case of Hamiltonian paths
(cf.\ \cite{BP,LM})

\vspace{1mm}

\noindent
{\bf Definition.}\ $\{ L_t \}$ is said to be {\it quasi-autonomous}
if there exist a lift $\{ \iota_t \}$ of $\{ L_t \}$ with its associated
function $h$ and the points $x_{+},x_{-} \in L$ such that
 \begin{eqnarray*}
  \max_{x \in L}h(t,x)=h(t,x_{+}), \ \
   \min_{x \in L}h(t,x)=h(t,x_{-})\\
 \end{eqnarray*}
and $\iota_{t}(x_{+}) \equiv \iota_{0}(x_{+}),\ \iota_{t}(x_{-}) \equiv \iota_{0}(x_{-})$
for all $t \in [0,1]$. 
And we call the lift $\{ \iota_t \}$ satisfying this property a
{\it quasi-autonomous lift}.

\vspace{1mm}

\noindent
{\bf Remark.}\ In \cite{Mi,Mi2} Milinkovi\'c studied geodesics in the space
$\mathcal{L}(T^{*}X,\omega,O_X)$ of Lagrangian submanifolds which are
Hamiltonian isotopic to the zero section $O_X$ of the cotangent bundle
$T^{*}X$ of a smooth closed manifold $X$.
In his terminology, the above condition is said to be
{\it strongly quasi-autonomous} (see \cite[Definition 10]{Mi2}).

\vspace{1mm}

The main theorem of this paper is the following.

\begin{thm}
Let $(M,\omega)$ be a symplectic manifold without boundary and
$L$ a closed Lagrangian submanifold of $M$.
Then, a regular exact path
$\{ L_t \}_{0 \leqq t \leqq 1} \subset \mathcal{L}(M,\omega,L)$
is $l$-critical if and only if it is quasi-autonomous.
\end{thm}

\noindent
{\bf Remark.}\ Take a lift $\iota_0$ of $L_0$.
If we choose a Hamiltonian path $\{ \psi_t \}_{0 \leqq t \leqq 1}$ such that
$L_t=\psi_{t}(L_0)$, then we have a lift $\{ \iota_t \}$ of the exact path
$\{ L_t \}$ such that $\iota_t=\psi_t \circ \iota_0$.
Its associated function is given by $h_t=H_t \circ \iota_t$.
By definition, Theorem 9 is equivalent to Theorem 1 in Section 1.

\vspace{1mm}

Let $\mathcal{V}_1$ be the tangent space at $0$ of the space $\mathcal{H}_1$
of all Hamiltonians generating loops in $\mathrm{Ham}^{c}(M,\omega)$, i.e.,
 \begin{eqnarray*}
  \mathcal{V}_1 := T_{0}(\mathcal{H}_1)
   =\Bigl\{ \frac{\partial F}{\partial s}(0) \Bigl|\ F \in \mathcal{V} \Bigl\}.
 \end{eqnarray*}
We quote a useful proposition (\cite[Proposition 12.2.B.]{Po}).
\begin{pro}[Polterovich]
The space $\mathcal{V}_1$ consists of all functions $G \in \mathcal{H}_0$
which satisfy
 \begin{eqnarray*}
  \int_{0}^{1} G(t,x)dt = 0
 \end{eqnarray*}
for all $x \in M$. \hfill \qed
\end{pro}

{\it Proof of Theorem 9.} Fix a lift $\{ \iota_t \}$ of $\{ L_t \}$ with
its associated function $h$.
We know from Lemma 8 that for any exact variation of $\{ L_t \}$
there exists $F \in \mathcal{V}$ such that the length function $l(s)$
has the form
 \begin{eqnarray}
  ||h-F(s)\circ\iota_t||
 \end{eqnarray}
and vice versa.
Moreover we know from the proof of Proposition 6 that $l(s)$ is approximated
up to the first order at $s=0$ by $||h-sF'(0)\circ\iota_t||$,
which is convex in $s$.
Thus with the Definitions after Proposition 6
it suffices to show the equivalence of the following two conditions:

\vspace{2mm}

$\bf (a)$ \ $\{ L_t \}$ is quasi-autonomous;

\vspace{1mm}

$\bf (b)$ \ For any fixed lift $\{ \iota_t \}$ of $\{ L_t \}$ and its
associated function $h$, the inequality

$||h-sG\circ\iota_t|| \geqq ||h||$ holds for any $G \in \mathcal{V}_1$
and $s \in \mathbb{R}$.

\vspace{2mm}

\noindent
$\bf (a) \Rightarrow (b):$ Since condition $\bf (b)$ is invariant under
the action of $\mathrm{Diff}(L)$, it suffices to show it for a specific lift.
Let $\{ \iota_t \}$ be a quasi-autonomous lift of $\{ L_t \}$ with
its associated function $h$.
For any $G \in \mathcal{V}_1$, we set
$u:=h-sG\circ\iota_t$. Then
 \begin{eqnarray*}
  ||u|| &=&
   \int_{0}^{1}\left( \max_{x \in L}u(t,x) - \min_{x \in L}u(t,x) \right)dt\\
  &\geqq& \int_{0}^{1}( u(t,x_+) - u(t,x_-) )dt\\
  &=& \int_{0}^{1}( h(t,x_+) - h(t,x_-) )dt
   -s \Bigl\{ \int_{0}^{1}G(t,\iota_{t}(x_+))dt
       - \int_{0}^{1}G(t,\iota_{t}(x_-))dt \Bigl\}.
 \end{eqnarray*}
Since $\{ L_t \}$ is quasi-autonomous, by Proposition 10, we have
 \begin{eqnarray*}
  \int_{0}^{1}G(t,\iota_{t}(x_{\pm}))dt
   = \int_{0}^{1}G(t,\iota_{0}(x_{\pm}))dt = 0.
 \end{eqnarray*}
Hence, we obtain $||u|| \geqq ||h||$, i.e., the condition $\bf (b)$.

To prove the converse, we must use a modified version of Akveld-Salamon's
extension of Hamiltonians.

Let $\{ L_t \}_{0 \leqq t \leqq 1}$ be an exact path of closed Lagrangian
submanifolds in a connected symplectic manifold $(M,\omega)$ without boundary.
Fix a lift $\{ \iota_t \}$ of $\{ L_t \}$ and denote its associated function
by $h_t$ and its push-forward via $\iota_t$ by $\tilde{h}_t$, i.e.,\
${\iota_t}_{*}h_t=\tilde{h}_t$.
Note that $\{ \tilde{h}_t \}$ depends only on $\{ L_t \}$,
not on the choice of the lift $\{ \iota_t \}$.
\begin{lem}
There exists a smooth extension 
$\{ H_t \}_{0 \leqq t \leqq 1} \subset C^{\infty}(M)$ of 
$\{ \tilde{h}_t \}_{0 \leqq t \leqq 1}$ such that 
 \begin{eqnarray*}
  \underset{M}{\mset} H_t = \underset{L_t}{\mset} \tilde{h}_t, \quad
  \underset{M}{\lset} H_t = \underset{L_t}{\lset} \tilde{h}_t,
 \end{eqnarray*}
for each $t \in [0,1]$.
Here, {\rm maxset} and {\rm minset} are defined as
 \begin{eqnarray*}
  & & \underset{M}{\mset}H_t:=\{ x \in M\ |\ H_t(x)=\max H_t \},\\
  & & \underset{M}{\lset}H_t:=\{ x \in M\ |\ H_t(x)=\min H_t \},\\
  & & \underset{L_t}{\mset}\tilde{h}_t
   :=\{ x \in L_t\ |\ \tilde{h}_t(x)=\max \tilde{h}_t \},\\
  & & \underset{L_t}{\lset}\tilde{h}_t
   :=\{ x \in L_t\ |\ \tilde{h}_t(x)=\min \tilde{h}_t \}
 \end{eqnarray*}
for each $t \in [0,1]$.
\end{lem}
{\it Proof.} We may assume that $\tilde{h}_t$ has positive maximum and
negative minimum if it is not constant, and that it is identically zero
if it is constant.

Since all $L_t$ are closed, we can choose a positive constant $\varepsilon$
such that $\mbox{exp}_x (v) $ gives a diffeomorphism between
the $\varepsilon$-neighbourhood of the zero section in the normal bundle
$N(L_t)$ and the $\varepsilon$-neighbourhood 
$U_{t,\varepsilon}$ of $L_t$ in $M$.
Take a nonnegative function $\alpha(t)\ (t \geqq 0)$
which is supported on $t \leqq \varepsilon$ and equals to $1$ 
on $t \leqq \varepsilon /2$.

Define $\{ H_t \}$ by
 \begin{eqnarray*}
  H_t(\mbox{exp}_x (v)) = \alpha (||v||^2) \tilde{h}_{t}(x) (1- || v ||^2),
   \quad x \in L_t, \quad v \in N_x(L_t)
 \end{eqnarray*}
on $U_{t,\varepsilon}$ and by $H_t=0$ outside $U_{t,\varepsilon}$.
This satisfies the required properties. \hfill \qed

\vspace{2mm}

Now we go back to the proof of Theorem 9.

\noindent
$\bf (b) \Rightarrow (a):$ Assume that $\bf (b)$ holds.
Setting $\tilde{h}_t:={\iota_t}_{*}h \in C^{\infty}(L_t)$,
this condition implies that
 \begin{eqnarray}
  \int_{0}^{1}
   \left\{ \max_{L_t}(\tilde{h}_t-G|_{L_t})
    - \min_{L_t}(\tilde{h}_t-G|_{L_t}) \right\}dt
  \geqq \int_{0}^{1}
    \left( \max_{L_t}\tilde{h}_t - \min_{L_t}\tilde{h}_t \right)dt
 \end{eqnarray}
for any $G \in \mathcal{V}_1$.
By Lemma 11, we can take a smooth extension of
$\{ H_t \}_{0 \leqq t \leqq 1} \subset C^{\infty}(M)$ of $\{ \tilde{h}_t \}$
such that $H_t \circ \iota_t=h_t$ and
 \begin{eqnarray}
  \underset{M}{\mset} H_t = \underset{L_t}{\mset} \tilde{h}_t, \quad
  \underset{M}{\lset} H_t = \underset{L_t}{\lset} \tilde{h}_t.
 \end{eqnarray}
We may assume that $H \in \mathcal{H}_0$.
Let us introduce a time-dependent smooth function $G$ on $M$ defined by
 \begin{eqnarray}
  G(t,y)=H(t,y)-\int_{0}^{1} H(t,y)dt.
 \end{eqnarray}
Then we have $G \in \mathcal{H}_0$ and
 \begin{eqnarray*}
  \int_{0}^{1} G(t,y)dt=0
 \end{eqnarray*}
for any $y \in M$.
Hence, Proposition 10 implies that $G \in \mathcal{V}_1$.

For this function $G$, by equations (12), (13) and (14), we have
 \begin{eqnarray*}
  ||H||
   &=& \int_{0}^{1}
    \left( \max_{L_t}\tilde{h}_t - \min_{L_t}\tilde{h}_t \right)dt\\
   &\leqq& \int_{0}^{1}
    \left\{ \max_{L_t}(\tilde{h}_t-G|_{L_t})
     - \min_{L_t}(\tilde{h}_t-G|_{L_t}) \right\}dt\\
   &=& \int_{0}^{1}
    \left\{ \max_{y \in L_t}(H(t,y)-G(t,y))
     - \min_{y \in L_t}(H(t,y)-G(t,y)) \right\}dt\\
   &=& \int_{0}^{1}
    \left\{ \max_{y \in L_t} \int_{0}^{1} H(t,y)dt
     - \min_{y \in L_t} \int_{0}^{1} H(t,y)dt \right\}dt\\
   &\leqq& \max_{y \in M} \int_{0}^{1} H(t,y)dt
     - \min_{y \in M} \int_{0}^{1} H(t,y)dt.
 \end{eqnarray*}
Hence we get
 \begin{eqnarray*}
  ||H||=\max_{y \in M} \int_{0}^{1} H(t,y)dt
     - \min_{y \in M} \int_{0}^{1} H(t,y)dt.
 \end{eqnarray*}
This equation yields that
 \begin{eqnarray*}
  \bigcap_{t \in [0,1]} \underset{M}{\mset}H_t \neq \phi,
   \quad
    \bigcap_{t \in [0,1]} \underset{M}{\lset}H_t \neq \phi.
 \end{eqnarray*}
(see e.g. \cite[Proposition 1.3.B.]{BP})

Since we took the extension as in Lemma 11,
the above condition is equivalent to
 \begin{eqnarray*}
  \bigcap_{t \in [0,1]} \underset{L_t}{\mset}\tilde{h}_t \neq \phi,
   \quad
    \bigcap_{t \in [0,1]} \underset{L_t}{\lset}\tilde{h}_t \neq \phi.
 \end{eqnarray*}
Conversely, suppose that 
$\bigcap_{t} \underset{L_t}{\mset}\tilde{h}_t$ and
$\bigcap_{t} \underset{L_t}{\lset}\tilde{h}_t$ are non-empty.
Take points $p_+ \in M$ in the former and $p_- \in M$ in the latter.
Since two curves
$(\iota_t)^{-1}(p_+)$ and
$(\iota_t)^{-1}(p_-)$ in $L$ are disjoint,
we can find a smooth family of diffeomorphisms
$\{ \phi_t \}_{0 \leqq t \leqq 1}$ of $L$
such that 
$(\iota_t \circ \phi_t )^{-1}(p_+)$ and 
$(\iota_t \circ \phi_t )^{-1}(p_-)$
are both fixed points in $L$, which we set $x_+$ and $x_-$, respectively.
We see that $\{ \iota_t \circ \phi_t \}$ is a lift of $\{ L_t \}$
with the associated function $h_t \circ \phi_t$,
which takes maximum at $x_+$ and minimum at $x_-$ in $L$.
Therefore, $\{ L_t \}$ is quasi-autonomous. \hfill \qed

\section{Applications and Examples}

\indent\indent
In this section we give some applications of our main theorem.

Theorem 1 says that all Lagrangian submanifolds $L_t$ in $M$ making
a $l$-critical exact path
$\{ L_t \}_{0 \leqq t \leqq 1} \subset \mathcal{L}(M,\omega,L)$
intersect each other.
Therefore we have

\begin{cor}
Let $L_0$ and $L_1$ be closed Lagrangian submanifolds in a symplectic manifold
$(M,\omega)$.
If $L_0 \cap L_1 = \phi$, then there exist no $l$-critical exact paths
connecting $L_0$ and $L_1$.
Hence there exist no length-minimizing geodesics between $L_0$ and $L_1$.
\end{cor}

We give an example of $l$-critical paths and pose a problem related
to it.
Let $L_0$ be the real projective space $\mathbb{R}P^n$ defined by
 \begin{eqnarray*}
  L_0:=\mathbb{R}P^n
   =\{ [z_0:\cdots:z_n]\ |\ z_0,\ldots,z_n \in \mathbb{R} \}.
 \end{eqnarray*}
in the complex projective
space $\mathbb{C}P^n$ endowed with the standard Fubini-Study symplectic
form $\omega_{FS}$.

Consider the Hamiltonian isotopy
 \begin{eqnarray*}
  \varphi_{t}^{H}([z_0:\cdots:z_n])
   =[z_0:e^{\pi i s t}z_1:\cdots:e^{\pi i s t}z_k:\cdots:z_n],
 \end{eqnarray*}
where $H_t=H: \mathbb{C}P^n \to \mathbb{R}$ is the autonomous Hamiltonian
given by
 \begin{eqnarray*}
  H([z_0:\cdots:z_n])
   =-\frac{s(|z_1|^2+\cdots+|z_k|^2)}{2(|z_0|^2+\cdots+|z_n|^2)} +const.
 \end{eqnarray*}
for $s \in [0,1]$ and $k \in \{ 1,2,\ldots,n \}$.

Put $L_t:=\varphi_{t}^{H}(L_0)$.
Then $\{ L_t \}_{0 \leqq t \leqq 1}$ is an exact path in 
$\mathcal{L}=\mathcal{L}(\mathbb{C}P^n,\omega_{FS},L_0)$.

A simple calculation yields
 \begin{eqnarray*}
  \underset{M}{\mset} H_t
   &=& \{ [z_0:\cdots:z_n]\ |\ z_0=\cdots=z_k=0 \} = \mathbb{C}P^{n-k},\\
  \underset{M}{\lset} H_t
   &=& \{ [z_0:\cdots:z_n]\ |\ z_0=0, z_{k+1}=\cdots=z_n=0 \} 
    = \mathbb{C}P^{k-1}.
 \end{eqnarray*}
Hence we have
 \begin{eqnarray*}
  \underset{M}{\mset} H_t = \mathbb{R}P^{n-k}, \quad
  \underset{M}{\lset} H_t = \mathbb{R}P^{k-1}
 \end{eqnarray*}
for all $t \in [0,1]$. Therefore,
 \begin{eqnarray*}
  \bigcap_{t} \underset{M}{\mset} H_t = \mathbb{R}P^{n-k}, \quad
  \bigcap_{t} \underset{M}{\lset} H_t = \mathbb{R}P^{k-1}
 \end{eqnarray*}
By theorem 1, we conclude that $\{ L_t \}$ is $l$-critical.

In the case where $s=1$, $\{ L_t \}_{0 \leqq t \leqq 1}$ is a loop in
$\mathcal{L}$.
Akveld-Salamon \cite{AS} proved that this loop is minimizing the Hofer
length over all periodic exact Lagrangian loops which are Hamiltonian
isotopic to $\{ L_t \} \subset \mathcal{L}$.

This result gives rise to the following natural question.

\vspace{1mm}

\noindent
{\bf Problem.}\ In the case where $s \in (0,1)$, is the above exact path
$\{ L_t \}$ a local minimizer of the Hofer's length function?,
or more strongly, is length-minimizing amongst all homotopic paths in
$\mathcal{L}$ with the same endpoints?

\section*{Acknowledgements}

We would like to thank Professors Paul Biran and Manabu Akaho for
valuable discussions.
The results of this paper were presented in the Symposium on Symplectic
Geometry and Topology in University of Tokyo (2006).
We thank Kaoru Ono, Martin Guest and Yukio Matsumoto for their warm
hospitality.

\begin{flushleft}
{\sc School of Engineering\\
Tokyo Denki University\\
Kanda-Nishiki-Cho, Chiyoda-Ku\\
Tokyo, 101-8457\\
Japan}

{\it e-mail} : {\tt hirie@im.dendai.ac.jp}
\end{flushleft}

\begin{flushleft}
{\sc College of Engineering\\
Nihon University\\
Koriyama\\
Fukushima, 963-8642\\
Japan}

{\it e-mail} : {\tt otofuji@ge.ce.nihon-u.ac.jp}
\end{flushleft}

\end{document}